\documentclass[12pt]{article}
\usepackage{amsmath,amssymb}
\usepackage{amsthm}
\usepackage{graphicx,tikz}
\usepackage{enumerate}
\usepackage{titlesec}
\usepackage{epstopdf}
\usepackage{caption}
\usepackage{xcolor}
\usepackage[sans]{dsfont}

\numberwithin{equation}{section}
\newtheorem{Thm}{Theorem}[section]
\newtheorem*{Thm*}{Theorem}
\newtheorem{Prop}[Thm]{Proposition}
\newtheorem{Lem}[Thm]{Lemma}

\newtheorem{Ques}[Thm]{Question}
\newtheorem{Cor}[Thm]{Corollary}
\newtheorem{Fact}[Thm]{Fact}

\newtheorem{Exa}[Thm]{Example}

\theoremstyle{remark}

\newtheorem{Rem}[Thm]{Remark}

\theoremstyle{definition}

\newtheorem{Def}[Thm]{Definition}

\newtheorem*{Def*}{Definition}

\numberwithin{equation}{section}

\newcommand{\m}[1]{\mathbb{ #1}}

     \def\ol{\overline}    
   
\def\al{\alpha}       \def\be{\beta}        \def\ga{\gamma}
\def\de{\delta}       \def\eps{\varepsilon}  
\def\th{\theta}       
\def\ka{\kappa}       \def\la{\lambda}      
\def\si{\sigma}                
\def\ph{\varphi}               
       \def\Ga{\Gamma}       
\def\La{\Lambda}

\theoremstyle{definition}

\theoremstyle{remark}

\newtheorem{Rmq}[Thm]{Remark}

\numberwithin{equation}{section}
\newfont{\goth}{eufm10 at 12pt}
\newfont{\gots}{eufm8 at 9pt}

\def\bt{\begin{Thm}}
\def\et{\end{Thm}}
\def\br{\begin{Rmq}}
\def\er{\end{Rmq}}

\def\bc{\begin{Cor}}
\def\ec{\end{Cor}}
\def\bp{\begin{Prop}}
\def\ep{\end{Prop}}
\def\bl{\begin{Lem}}
\def\el{\end{Lem}}
\def\bd{\begin{Def}}
\def\ed{\end{Def}}
\def\bq{\begin{quotation}}
\def\eq{\end{quotation}}
\def\bfa{\begin{Fact}}
\def\efa{\end{Fact}}
\def\bexa{\begin{Exa}}
\def\eexa{\end{Exa}}

\def\ra{\rightarrow}
\def\vs{\vspace{1em}}

\begin{document}
\title{
Convolution and square 
in abelian groups I
}
\author{Yves Benoist
}
\date{}

\maketitle

\begin{abstract}
\noindent 
We prove that the  functional equation $f\!\star\! f (2\,t) =\la f(t)^2$,
for $t$ in $\m Z/d\m Z$ with $d$ odd, admits a non-zero 
solution $f$ if $\la=\sqrt{a}\!+\! i\sqrt{b}$ with $a,b$ positive integers such that $a\!+\! b\!=\!d$ and 
$a\!\equiv\frac{(d+1)^2}{4}$ mod $4$.
The proof involves theta functions on elliptic curves with complex multiplication.
\end{abstract}

\renewcommand{\thefootnote}{\fnsymbol{footnote}} 
\renewcommand{\thefootnote}{\arabic{footnote}}

\section{The functional equation}

It is well known that elliptic curves, theta functions and modular forms are 
very useful tools both in Algebraic Geometry and in Number Theory -see \cite{Farkas08}.
In this paper, we focus on an elementary open problem 
raised by numerical experiments.  The main surprise 
in this text is that the intriguing output of these numerical experiments
cannot be understood without these tools.

\subsection{Introduction}

We will
deal with a finite abelian group $G$ of odd order $d$, 
which, most of the time, will be the cyclic group $G=\m Z/d\m Z$, and with the functional equation 
\begin{equation}
\label{eqntante} 
f*f\,(2t) =\la \, f^2(t)
 \;\;\; \mbox{\rm for all $t$ in $G$,}
\end{equation}
where the unknown  
is a non-zero function $f: G\ra \m C$
and where $\la\in \m C$ is a parameter.
This equation expresses 
a proportionality condition between the ``convolution square'' of $f$ and 
its ``multiplication square''. 
A non-zero solution $f$ of this functional equation \eqref{eqntante} will be called a ``$\la$-critical
function on $G$'' or, in short, a ``$\la$-critical
function'',
and a value $\la$ for which such a function $f$ 
exists will be called a ``critical value on $G$'', or a ``$d$-critical value''
when $G=\m Z/d\m Z$. 
Note that Equation \eqref{eqntante} has been chosen so that it is invariant by translation on the variable $t$. 
This equation \eqref{eqntante} 
can be rewritten as 
\begin{equation}
\label{eqntante2} \textstyle
\sum\limits_{\ell\in G}f(k\!+\!\ell)\,f(k\!-\!\ell) =\la \, f(k)^2
 \;\;\; \mbox{\rm for all $k$ in $G$.}
\end{equation}

The aim of this text is to point out 
the interest of this functional equation
by gathering unexpected results and questions based on numerical experiments
and by relating this apparently naive functional equation to elliptic curves with complex multiplication. 

Indeed, our main result, Theorem \ref{thmrapirb}, gives explicit $d$-critical values.
Eventhough the statement of this theorem
is very short and purely elementary, surprisingly, 
our construction relies on the
Jacobi theta functions $\theta(z,\tau)$ for special values of the parameter $\tau$,
and the key point in the proof relies on  modularity properties due to Hecke
(Lemma \ref{lemtrafor}) of these theta functions.
\vs

I tried to keep this text as elementary and concrete as possible. In a second more technical paper \cite{CSAGII}, 
I will extend this construction
of critical values to all finite abelian groups $G$
by using  
the Riemann theta functions on higher dimensional abelian varieties, 
and their modularity properties
as functions on the Siegel upper half-space -- see Theorem \ref{thmmaivap}.
\vs 

I thank G\'{e}rard Laumon, Samuel Leli\`evre and Emmanuel Ullmo  for enlighting discussions
on this project.

\subsection{Comments}
\hspace*{1.5em}$\star$ 
Since $d$ is odd, the value $\la=0$ is not $d$-critical. Indeed
the only function $f$ for which $f*f=0$ is $f=0$.

$\star$ By analogy, one might look at Equation \eqref{eqntante} on  locally compact abelian groups $G$. When the group is $G=\m R$
or $G=\m Z$, the only $L^2$-solutions that I know are gaussian functions $f(t)=e^{at^2+bt+c}$ where $a,b,c \in \m C$,  $Re(a)<0$, 
together with, when $G=\m Z$, their restrictions to subgroups. 
When the group is $G=\m R/\m Z$, the only $L^2$-solutions that I know are constant functions.

\subsection{Special critical values}
\label{secspecri}
\noindent 

In this section, we list a few $d$-critical values that are easy to find.
We call them special.

$\star$ If we choose $f(0)=1$ and $f(k)=\al$ to be constant for $k$ in  
$G\smallsetminus\{0\}$, we find four critical values:
we find $\la = 1$ when $\al=0$,  
we find $\la=d$ when $\al=1$, and
we find $\la=\frac{d-3+\eps\sqrt{D}}{2}$ when $\al=\frac{1-d-\eps\sqrt{D}}{2(d-1)}$,
with $\eps=\pm 1$ and $D=(d\!-\!1)(d\!-\!9)$. Note that these $d$-critical values are real as soon as $d\geq 9$.

$\star$ If we choose $f$ to be a gaussian function
$f(k):=\eta^{-k^2}$ with $\eta:=-e^{i\pi/d}$, we find 
the critical value $\la=\sqrt{d}$ when $d\equiv 1\,{\rm mod}\, 4$ and $\la= i\sqrt{d}$ when $d\equiv 3$ mod $4$.
Moreover, its opposite $-\la$ is also often a critical value. 
This is the case when $d$ is not a square.
But this is not always the case, for instance, when $d=9$, the value $-3$ is not
$d$-critical.

\subsection{Induced critical values}
\noindent

In this section we explain how to construct $d$-critical values when $d$ is a composite number starting from critical values for the factors of $d$. 
The method works for any abelian groups.

$\star$ Let $G_1\subset G$ be finite abelian groups and $\la_1$
be a critical value on $G_1$ with $\la_1$-critical function $f_1$.
Then $\la_1$ is also a critical value on $G$ 
with $\la$-critical function $f:=f_1{\bf 1}_{G_1}$. 

$\star$ Let $G_1\subset G$ be finite abelian groups, $d_1$ be the order of $G_1$ and $\la_2$ be a critical value on the quotient $G/G_1$ 
with $\la_2$-critical functions $f_2$.
Then $\la:=d_1\la_2$ is a critical value on $G$ 
with $\la$-critical function $f:=f_2\circ \pi$
where $\pi:G\ra G/G_1$ is the projection.

$\star$ Let $G=G_1\times G_2$ be the product 
of two finite abelian groups,
let $\la_1$
be a critical value on $G_1$ with $\la_1$-critical function $f_1$
and $\la_2$
be a critical value on $G_2$ with $\la_2$-critical function $f_2$.
Then the product 
$\la:=\la_1\la_2$ is a critical value on $G$ with $\la$-critical function 
$f:=f_1\otimes f_2$.
\vs 

The most interesting critical values will be those that are not special and that are not obtained
by these ``induction'' methods.

\subsection{Numerical experiments}
\label{secnumexp}

\noindent 
The following lists of $d$-critical values rely on numerical experiments using the Buchberger's
algorithm for computing Groebner basis (see \cite[Chap. 2]{CLO92}).
We denote by $b_d$ the size of the list.
\vs

For $3\leq d\leq 9$ the complete lists of $d$-critical values are: 
\vs 

\noindent 
$\star$ \fbox{$d=3$} $b_3=4$.\; $\la=1$, $3$, and $\pm i\sqrt{3}$.

\noindent 
$\star$  \fbox{$d=5$} $b_5=6$.\; $\la=1$, $5$, $\pm \sqrt{5}$, and $1\pm 2i$.

\noindent 
$\star$  \fbox{$d=7$} $b_7=8$.\; $\la=1$, $7$, $\pm i \sqrt{7}$, and $\pm 2\pm i\sqrt{3}$.

\noindent 
$\star$ \fbox{$d=9$} $b_9\!= \!15$.\; $\la=1$, $9$, $3$, $\pm i\sqrt{3}$, $\pm 3i\sqrt{3}$, $\pm 1\pm 2i\sqrt{2}$, and $\pm \sqrt{5}\pm 2i$.
\vs 

\noindent For $d=11$ and $ 13$, the  lists below of $d$-critical values are still probably complete. 
\vs 

\noindent $\star$ \fbox{$d=11$} $b_{11}\!=\!20  $.\; $\la =1$ , $11$, $4\pm \sqrt{5}$, $\pm i\sqrt{11}$,\\
\hspace*{1em} $2\pm i\sqrt{7}$, $\pm 2\sqrt{2}\pm i\sqrt{3}$, and $\pm (1\!+\!\eps\sqrt{5})\pm i\sqrt{5\!-\!2\eps \sqrt{5}}$ with $\eps=\pm 1$.
\vs 

\noindent $\star$ \fbox{$d=13$} $b_{13}\!=\! 18$.\; $\la =1$ , $13$, $5\pm 2\sqrt{3}$, $\pm i\sqrt{13}$,\\
\hspace*{1em} $\pm 1\pm 3i\sqrt{2}$, $\pm \sqrt{5}\pm 2i\sqrt{2}$, and $\pm 3\pm2i$.
\vs 

\noindent For $d=15$ and $17$, here are  just a few $d$-critical values. These values were obtained by the guess and check method. Looking at the critical values for $d\leq 13$ one can guess a few critical value for $d=15$ and $17$. The key point of the method 
is that  the Buchberger's algorithm is much faster when  checking if a given $\la$ is $d$-critical than when finding all the $d$-critical values.
\vs

\noindent $\star$ \fbox{$d=15$} $b_{15}\!=\! 60$.\;
$\la=$ product of a $3$-critical and a $5$-critical value, and\\
\hspace*{1em} $\la=-3$, $-5$, $6\pm \sqrt{21}$,
$\pm 2\pm i\sqrt{11}$, $\pm 2\sqrt{2}\pm i\sqrt{7}$, $\pm 2\sqrt{3}\pm i\sqrt{3}$, and\\
\hspace*{1em}  $\pm \!2\sqrt{2\!-\!\eps\sqrt{3}} \pm (2\!+\eps\!\sqrt{3})i$ and $1\!+\!\eps\sqrt{5}\pm i\sqrt{9\!-\!2\eps \sqrt{5}}$ with $\eps=\pm 1$, and\\ 
\hspace*{1em} $\pm(\sqrt{3}\pm i\sqrt{2})(\sqrt{2}\pm i)$.
\vs 

\noindent $\star$ \fbox{$d=17$} $b_{17}\!=\! 28$.\; $\la =1$ , $17$, $7\pm 4\sqrt{2}$, $\pm i\sqrt{17}$,\\
\hspace*{1em} $\pm 1\pm 4i$, $\pm \sqrt{5}\pm 2i\sqrt{3}$, $3\pm 2i\sqrt{2}$, $\pm \sqrt{13}\pm 2i$, and\\
\hspace*{1em} $\pm(1+2\eps\sqrt{2})\pm 2i\sqrt{2-\eps\sqrt{2}}$ with $\eps=\pm 1$.\vs 

\noindent
For our numerical experiments we used SageMath and Maple softwares. 

\section{ Critical values}

One of the motivations of Proposition \ref{procrival} and  Theorem \ref{thmrapirb}
below is to explain some of 
the intriguing patterns that occur in these experimental  lists of critical values.

\subsection{Properties of critical values} 

We first begin by a few properties of the critical values, 
that are valid on any finite abelian group. 

\bp
\label{procrival} Let $G$ be a finite abelian group of odd order $d$, and $\la$ a critical value on $G$, then:\\ 
$(i)$ all the Galois conjugate of $\la$ are also critical values on $G$,\\
$(ii)$ one has $|\la|\leq d$ with equality if and only if $\la=d$,\\ 
$(iii)$ the ratio $d/\la$ is also a critical value on $G$,\\
$(iv)$ The ratio $\frac{\la -1}{2}$ is an algebraic integer.
\ep

\bc 
\label{corcrival} There exist only finitely many critical values on $G$.
\ec 

\begin{proof}[Proof of Corollary \ref{corcrival}]
Since it is obtained by elimination, the set of critical values on $G$ is either finite or its complement in $\m C$ is finite. Since, by $(ii)$ it is bounded, it must be finite.
\end{proof}

\begin{proof}[Proof of Proposition \ref{procrival}] 
$(i)$ Equations \eqref{eqntante2} have rational coefficients.

$(ii)$ This follows from Cauchy-Schwarz inequality. Indeed, setting\\ 
$\|f\|_\infty=\max\limits_{k\in G}|f(k)|$ and $\|f\|_2=(\sum_k|f(k)|^2)^{\frac12}$,
one has
$$
|\la|\|f\|^2_\infty=\| f\star f\|_\infty\leq \| f\|_2^2
\leq d\, \| f\|_\infty^2.
$$
Hence $|\la|\leq d$. In case we have equality the function $f$ must have constant modulus,
and must satisfy
$f(k\!+\!\ell)f(k\!-\!\ell)=f(k)^2$, for all $k$, $\ell$. Hence $f$ is proportional to a character
and one has $\la=d$.

$(iii)$ If $f$ is a $\la$-critical function on $G$, then its Fourier transform $\widehat{f}$, which is given by, for every character $\chi:G\ra \m C^*$, 
$$
\textstyle
\widehat{f}(\chi)= \sum\limits_{x\in G}f(x)\chi(x),
$$ 
is a $d/\la$-critical function on the dual group $\widehat{G}$
which is isomorphic to $G$.

$(iv)$ Let $G_+$ be a subset of $G$ of cardinality $\frac{d-1}{2}$
such that for each non-zero $\ell\in G$ either $\ell$ or $-\ell$ is in $G_+$.  The equations \eqref{eqntante2} can be rewritten as
\begin{equation}
\label{eqntante3} 
\tfrac{\la -1}{2} \, f(k)^2 =
\textstyle\sum\limits_{\ell\in G_+}f(k\!+\!\ell)\,f(k\!-\!\ell) 
 \;\;\;\;\;\; \mbox{\rm for all $k$ in $G$}
\end{equation}
Let $K$ be the subfield of $\m C$ generated by the coefficients $f(k)$.
To prove that  $\la':=\tfrac{\la -1}{2}$ is an algebraic integer, it is enough to check that, for all non-archimedean absolute value
$|.|_v$ on $K$, one has $|\la'|_v\leq 1$.\\ 
We set $\| f\|_v:= \max\limits_{\ell\in G} |f(\ell)|_v$, we choose $k$ such that 
$\|f\|_v=|f(k)|_v$, and we compute
\begin{eqnarray*}
|\la'|_v \|f\|_v^2=|\la'f(k)^2|_v&=&|\textstyle\sum\limits_{\ell\in G_+} f(k\!+\!\ell)f(k\!-\!\ell)|_v\\
&\leq& \max\limits_{\ell\in G}\, |f(k\!+\!\ell)|_v|f(k\!-\!\ell)|_v
\;\leq\; \|f\|_v^2.
\end{eqnarray*}
This proves that $|\la'|_v\leq 1$ as required.
\end{proof}

\subsection{Construction of critical values} 

From now on, $G$ will be the cyclic group $\m Z/d\m Z$. 
It is not clear from the definition that there does exist 
$d$-critical values that are non-induced and non-special. 
The following theorem tells us that this is always the case
for $d\geq 5$.

\bt
\label{thmrapirb}
Let $a$,$b$ be positive integers with 
$a\!+\!b\!=\!d$ and  
$a\!\equiv\!\frac{(d+1)^2}{4}\;{\rm  mod}\; 4$.
Then the complex number $\la:=\sqrt{a}+i\sqrt{b}$
is a $d$-critical value.
\et

\begin{Rem} The congruence assumption in Theorem \ref{thmrapirb} is equivalent to
\begin{eqnarray}
\label{eqnabmabm}
a-b\equiv 1 \;\mbox{\rm mod}\; 4 
&{\rm and}&
ab \equiv 0\;\mbox{\rm mod}\; 4.
\end{eqnarray}
A more concrete way to state Theorem \ref{thmrapirb} is:

For $d\equiv 1$ mod $4$, the following values are $d$-critical:\\
$\sqrt{d}$ ,\; $\sqrt{d\!-\!4}\!+\!2i$ ,\; $\sqrt{d\!-\!8}\!+\! 2i\sqrt{2}$ ,\; $\sqrt{d\!-\!12}\!+\!2i\sqrt{3}$ , ...
 
 For $d\equiv 3$ mod $4$, the following values are $d$-critical:\\
$i\sqrt{d}$ ,\; $2\!+\!i\sqrt{d\!-\!4}$ ,\;  $2\sqrt{2}\!+\!i\sqrt{d\!-\!8} $ ,\;  $2\sqrt{3}+i\sqrt{d\!-\!12}$ , ...
\end{Rem} 

More precisely, we will see that, surprisingly, for these values $\la$,
the set of $\la$-critical functions has positive dimension. Indeed, 
we will construct a one-parameter family of 
$\la$-critical functions using a suitable Jacobi theta function.

Before that we discuss the above congruence condition on $a$.

\bl
\label{lemalgint}
Let $a$,$b$ be positive integers with 
$a\!+\!b\!=\!d$ and let  $\la:=\sqrt{a}+i\sqrt{b}$.
The number $\frac{\la-1}{2}$ is an algebraic integer if and only if
$a\!\equiv\!\frac{(d+1)^2}{4}\;{\rm  mod}\; 4$.
\el

In particular, by Proposition \ref{procrival}.$iv$, when 
$a\!\not\equiv\!\frac{(d+1)^2}{4}\;{\rm  mod}\; 4$,
the complex number $\la=\sqrt{a}+i\sqrt{b}$
can not be a $d$-critical value.
\vs 

\begin{Rem}
Note that, for any algebraic number $\la$, one has the equivalence:
\begin{eqnarray}
\nu\!:=\!\tfrac{\la-1}{2}\; \mbox{\rm is an algebraic integer}
\;\Longleftrightarrow\; 
\nu'\!:=\!\tfrac{\la^2-1}{4}\; \mbox{\rm is an algebraic integer}.&&
\;\;\mbox{ }
\end{eqnarray}
Indeed, theses two elements $\nu$ and 
$\nu'$ are related by the equation
$\nu^2+\nu=\nu'$.
\end{Rem}

\begin{proof}[Proof of Lemma \ref{lemalgint}] 
The number $\nu'\!:=\!\tfrac{\la^2-1}{4}$ is equal to
$\nu'=\tfrac{a-b-1}{4}+i\tfrac{\sqrt{ab}}{2}$.
It is an algebraic integer if and only if 
$a-b\equiv 1 \;\mbox{\rm mod}\; 4$ 
\; {\rm and}\;
$ab \equiv 0\;\mbox{\rm mod}\; 4$.
As  seen in \eqref{eqnabmabm}, this condition is equivalent to 
$a\!\equiv\!\frac{(d+1)^2}{4}\;{\rm  mod}\; 4$.
\end{proof} 

\bc
\label{corrapirb2}
Let $p$,$q$ be positive integers with $p$ odd and $q$ even and let $d:=p^2+q^2$.
Then the complex number $\la:=p+iq$
is a $d$-critical value.
\ec

\begin{proof} Condition \eqref{eqnabmabm} is true:
$p^2-q^2\equiv 1$ mod $4$\; and\; $p^2q^2 \equiv 0$ mod $4$.
\end{proof}

\subsection{More numerical experiments}

A reasonable aim in this topic would be to give, 
for each $d$ the list of the $d$-critical values $\la$,
and for each $\la$ the description of the projective algebraic variety 
given by the $\la$-critical functions (number of connected components,
their dimension,...).

Here are a few less ambitious questions 
supported by numerical experiments
that suggest that some hidden structure has to be understood. 
\vs 
 
The first question deals with the properties of the $\la$-critical functions
when $\la$ belongs to a real or imaginary quadratic number field.

\begin{Ques}
\label{concrisym} 
Let $d$ be an odd integer and $\la$ a $d$-critical value which is quadratic.
Does  there exist an even $\la$-critical function?
\end{Ques}

That is a $\la$-critical function $f$ such that $f(-k)=f(k)$ for all $k$.

We checked, using numerical experiments, that this is true for $d\leq 11$.
Unfortunately this is not true when $\la$ is not quadratic as 
$2\sqrt{2}-1+2i\sqrt{2\!+\!\sqrt{2}}$.\vs 

The second question deals with the properties of the $d$-critical values.

\begin{Ques}
\label{condslgal} Let $d$ be an odd prime and $\la$ a $d$-critical value.\\ 
Are $\la$ and $d/\la$ Galois conjugate, except for $\la=1$ and $\la=d$?
\end{Ques}

More generally, when $d$ is not prime,  one might still expect a similar question to be true for {\it non-induced} $d$-critical values.

Note that we do not expect all the Galois conjugates in $\m C$ of the 
non-real critical values $\la$ to
have absolute value equal to  $\sqrt{d}$. 
We computed, using numerical experiments,  an example of $d$-critical value $\la$ of degree $8$ over $\m Q$ 
with two real and six non-real Galois conjugates in $\m C$. 
\vs 

The third question deals with critical values that are real quadratic.

\begin{Ques}
\label{conquacri} 
Let $d$ be an odd prime
and $\la$  a $d$-critical value 
which is a real quadratic number. Is such a $\la$   special?
\end{Ques}

We recall that special means that $\la=1$, $d$, $\pm \sqrt{d}$, or $\frac{d-3\pm\sqrt{(d-1)(d-9)}}{2}$.

Note that the $d$-critical values which are quadratic over $\m Q$ and non-real
can be described up to sign thanks to 
Proposition \ref{procrival} and Theorem \ref{thmrapirb}. 

We checked, using numerical experiments, that this  is true for $d\leq 11$.
\vs 

The last question deals with critical values that are quadratic over $\m Q$ but are non real. More precisely it deals with the sign of $Re(\la)$ in Theorem \ref{thmrapirb}.

\begin{Ques}
\label{consgncri} Let $d$ be an odd integer and  $a$,$b$ be positive integers with 
$a\!-\! b\!\equiv\! \frac{d^2+1}{2}\,{\rm mod}\, 8$ and $a\!+\!b\!=\!d$.
If $d\not\equiv 2$ mod $3$, 
is the number $\la:=-\sqrt{a}+i\sqrt{b}$
a $d$-critical value?
\end{Ques}

Moreover, for each $d\equiv 2$ mod $3$, 
one still expects the answer to this question to be true 
with at most one exception.
Recall that the interesting case is 
when the critical value $\la$ is quadratic over $\m Q$ i.e. when the integer $a$ is a square (see Remark \ref{remmla}).

We checked, using numerical experiments,  that
$-1+2i$ and $-2+i\sqrt{7}$ are not  critical values
and that, for $d\leq 23$, 
the only other possible exceptions are $-3+2i\sqrt{2}$ 
and $-4+i\sqrt{7}$. 
\vs 

There is another similar question: 
let $d=a^2$ be the square of an odd integer $a\geq 3$. 
We know that this integer $a$ is a $d$-critical value. 
But when is its opposite $-a$ also $d$-critical?  I checked, using numerical experiments, that $-3$ is not $9$-critical  but that $-5$ is $25$-critical.

\section{Theta functions and elliptic curves}
Our aim now is to prove Theorem \ref{thmrapirb}.

\subsection{Main result} 

We recall the definition of the Jacobi theta function:
$$\textstyle
\theta_\tau (z)=\theta(z,\tau):=\sum\limits_{m\in \m Z}e^{i\pi\tau m^2}e^{2i\pi mz},
\;\;\mbox{\rm for $z\in \m C$ and $\tau\in \m H$,}
$$
where $\m H$ is the upper half plane 
$\m H=\{ \tau\in \m C\mid Im(\tau)>0\}$. 
This function is $1$-periodic: $\theta_\tau(z+1)=\theta_\tau(z)$.
We can now explain our construction of $\la$-critical functions
on $\m Z/d\m Z$.

\bd
\label{defrapirb}
We will say that the function $\theta_{\tau}$ is $(\la,d)$-critical 
if, for all $z$ in $\m C$, the function
$
f_{z,\tau}:\ell\mapsto \theta(z\!+\!\ell/d,\tau)
$ 
is $\la$-critical on $\m Z/d\m Z$.
\ed

This means that, for all $z$ in $\m C$,
$$\textstyle
\sum\limits_{\ell\in\m Z/d\m Z} \theta(z+\ell/d,\tau)\,\theta(z-\ell/d,\tau)
\;=\;\la\,\theta(z,\tau)^2.
$$

Theorem \ref{thmrapirb} is a special case of the following Proposition \ref{prorapirb}.$(a)$. 
The whole Proposition \ref{prorapirb} tells us more. It tells us exactly for which parameters  $d$, $\la$, $\tau$,  the function $\theta_\tau$ is $(\la,d)$-critical. 

\bp
\label{prorapirb}
Let $a$,$b$ be positive integers with 
$a\!\equiv\!\frac{(d+1)^2}{4}\;{\rm  mod}\; 4$
and $a\!+\!b\!=\!d$.
Set $\la_0:=\sqrt{a}+i\sqrt{b}$ and
\begin{eqnarray}
\label{eqntau0}
\tau_0&:=&\tfrac{1}{4d^2}(a-b-d^2+2i\sqrt{ab}).
\end{eqnarray}
$(a)$ The function $\theta_{\tau_0}$ is $(\la_0,d)$-critical.\\
$(b)$ Conversely, let $\tau\in \m H$ such that the function $\theta_\tau$
is $(\la,d)$-critical, 
then one has $\la=\pm\sqrt{a}\pm i\sqrt{b}$ with $a$, $b$ as above.\\
$(c)$ The function $\theta_\tau$ is 
$(\la,d)$-critical for $\la=\pm\la_0$
if and only if $\tau=(k+\tau_0)/p$ with $\tau_0$ as above, $k\in \m Z$ and $p>0$ a 
divisor of the integer $N_k:=d^2\, |k\!+\! \tau_0|^2$.\\ 
$(d)$ The above sign $\eps=\pm$ is given by the Jacobi symbol $\eps=(\!\frac{p}{4k-1}\!)$. 
\ep

Remember that the Jacobi symbol $(\!\frac{\ga_0}{\de_0}\!)=\pm 1$ is defined 
for two relatively prime integers $\ga_0$ and $\de_0$ with $\de_0$ odd, and that,
by convention, when $\de_0$ is negative, it is given by $(\!\frac{\ga_0}{\de_0}\!)=(\!\frac{\ga_0}{-\de_0}\!)$.

Note also that, in view of Point $(b)$, the assumption $\la=\pm \la_0$ in Point $(c)$ is not restrictive since one has the equivalence: 
$$\mbox{$\theta_\tau$ is $(\la,d)$-critical\; $\Longleftrightarrow$ \;
$\theta_{-\ol{\tau}}$ is $(\ol{\la},d)$-critical.}
$$ 

\br The parameter $\tau_0$ will be called the {\it fundamental parameter} and the parameters 
$\tau_{k,p}:= (k+\tau_0)/p$ 
the {\it associated parameters}.
These parameters $\tau_{k,p}$ and the integers $N_k$ can also be given 
by the simple formulas with $m_0:=\frac{a-b-d^2}{4}$ 
and $N_0:=\frac{(d+1)^2-4a}{16}$:
$$
\fbox{$\tau_{k,p}=\frac{1}{d^2p}(d^2k+m_0+i\sqrt{ab/4})$}
\;\; {\rm  where}\;\; 
\fbox{$p\, | \, N_k:=d^2k^2+2m_0k+N_0$}.
$$ 

\er

\bexa To be very concrete, we give below 
the list of all values 
$\tau=\tau_{k,p}$ for which $\theta_\tau$ is $(\pm\la_0,d)$-critical  with $k$ in $\m Z$
and $p$ divisor of $N_k$, when $d\leq 9$.\vs

\noindent
$\star$  $d=5$, $\la_0=1+2i$\;\; : 
$\tau_{k,p}=\tfrac{1}{25 p}(25k - 7+i)$\;\; \;\, 
where\; 
$p\,|\, 25k^2\!-\! 14k\!+\! 2$.\\
$\star$  $d=7$, $\la_0=2+i\sqrt{3}$: 
$\tau_{k,p}=\tfrac{1}{49 p}(49k\!-\! 12+i\sqrt{3})$ where\; $p\,|\, 49k^2\!-\! 24k\!+\! 3$.\\
$\star$  $d=9$, $\la_0=1\!+\! 2i\sqrt{2}$: 
$\tau_{k,p}=\tfrac{1}{81 p}(81k\!-\! 22+i\sqrt{2})$ where\; $p\,|\, 81k^2\!-\! 44k\!+\! 6$.\\
$\star$  $d=9$, $\la_0=\sqrt{5}+2i$: 
$\tau_{k,p}=\tfrac{1}{81 p}(81k\!-\! 20+i\sqrt{5})$ where\; $p\,|\, 81k^2\!-\! 40k\!+\! 5$.
\eexa

Given $d$ and $\la_0$, we have  seen that it is  always possible to 
choose $k$ and $p$ such that  $\la=\eps\la_0$ with  sign $\eps=(\!\frac{p}{4k-1}\!)$ equal to $+1$: we just choose $p=1$.

On the other hand, given $d$ and $\la_0$, it is sometimes possible to choose 
$k$ and $p$ such that $\la=\eps\la_0$ with  sign $\eps=(\!\frac{p}{4k-1}\!)$ equal to  $-1$. For instance,\\  
$\star$ when $\la= -2-i\sqrt{3}$  with
$\tau=\tfrac{1}{98}(37+i\sqrt{3})$, or \\
$\star$ when $\la=-\sqrt{5}-2i$ with
$\tau=\tfrac{1}{162}(61+i\sqrt{5})$.\\
We mention the following corollary  that tells us exactly when this  is possible.

\bc 
\label{corrapirb}
Let $a$,$b$ be positive integers with 
$a\!\equiv\!\frac{(d+1)^2}{4}\;{\rm  mod}\; 4$
and $a\!+\!b\!=\!d$.
Set $\la:=-\sqrt{a}-i\sqrt{b}$. 
There exists $\tau$ in $\m H$ whose function $\theta_\tau$
is $(\la,d)$-critical if and only if either
$a$ is not a square in $\m Z$ or $-b$ is not a square in $\m Z/2a\m Z$.
\ec

Concretely, here are a few special cases of our criteria for such a function $\theta_\tau$ to exist, in which $a$ is a square in $\m Z$ and $b$ is either  $4\ell$ or $4\ell-1$ with $\ell$  a positive integer:\\
$\star$ When $\la=-1-2i\sqrt{\ell}$ : never.\\
$\star$ When $\la=-3-2i\sqrt{\ell}$ : $\ell\equiv 1$ mod $3$.\\
$\star$ When $\la=-5-2i\sqrt{\ell}$ : $\ell\equiv 2$ or $3$ mod $5$.\\
$\star$ When $\la=-2-i\sqrt{4\ell-1}$: $\ell\equiv 1$ mod $2$.\\
$\star$ When $\la=-4-i\sqrt{4\ell-1}$: $\ell\equiv 1$ mod $2$.\\
$\star$ When $\la=-6-i\sqrt{4\ell-1}$: $\ell\equiv 1,2,3$ or $5$ mod $6$.

The proof of this corollary that relies on quadratic reciprocity is left to the reader.
We will not use it below. The rest of the paper is devoted to the proof of Proposition
\ref{prorapirb}.

\subsection{Preliminary formulas} 

The proof of Proposition \ref{prorapirb}.$(a)$ relies on three classical formulas for the theta functions,
the ``addition formula'', the ``isogeny formula'',  the ``transformation formula''. We will only need special cases of these 
formulas that we state below.

We need to introduce the theta 
functions\footnote{\noindent 
With other ''classical'' notations
for theta functions as in \cite{BiLa04}, one has the equalities,\\
$\theta_{[0]}(z,\tau)=\theta_{0,0}(2z,2\tau)=
\theta
\mbox{\scriptsize 
$\left[\!\!\!
\begin{array}{c} 0 \\0\end{array}\!\!\!
\right]$} 
(2z,2\tau)$ 
and
$\theta_{[1]}(z,\tau)=\theta_{1,0}(2z,2\tau)=
\theta
\mbox{\scriptsize 
  $\left[\!\!\!\!
   \begin{array}{c} 1/2\\0\end{array}\!\!\!\!
   \right]$} 
(2z,2\tau).$
}
\begin{eqnarray*}
\theta_{[0]}(z )=\theta_{[0]}(z,\tau)
&:=&\textstyle
\sum\limits_{m\,\rm even}e^{i\pi\frac{\tau}{2} m^2}e^{2i\pi mz}
\\
\theta_{[1]}(z )=\theta_{[1]}(z,\tau)
&:=&\textstyle
\sum\limits_{m\,\rm odd}e^{i\pi\frac{\tau}{2} m^2}e^{2i\pi mz}.
\end{eqnarray*}

Note that one has the equalities:
\begin{eqnarray}
\label{eqnt0t1t2}\textstyle
\theta_{[0]}(z,\tau)=\theta(2z,2\tau)
&\;{\rm and}\;&
\theta_{[0]}(z,\tau)+\theta_{[1]}(z,\tau)=\theta(z,\tau/2).
\end{eqnarray}

Here is the first formula that we need.

\bl 
\label{lemaddfor} {\bf Addition formula} For all $z,w$ in $\m C$, $\tau\in \m H$, one has
\begin{eqnarray}
\label{eqnaddfor}
\theta(z+w,\tau)\theta(z-w,\tau)&=&
\theta_{[0]}(w,\tau )\theta_{[0]}(z,\tau )+
\theta_{[1]}(w,\tau )\theta_{[1]}(z,\tau ).
\end{eqnarray}
\el 

\begin{proof} 
Just write the left-hand side as a double sum over $m$, $n$ in $\m Z$
and split this double sum according to the parity of $m\!-\! n$.
\end{proof}

Here is the second  formula which is simple but useful.

\bl 
\label{lemisofor} {\bf Isogeny formula} For $\tau\in \m H$, $d$ odd positive integer, one has
\begin{eqnarray*}
\label{eqnisofor}\textstyle
\sum\limits_{\ell \in \m Z/d\m Z}\theta_{[0]}(\ell /d,\tau)=d\,\theta_{[0]}(0,d^2\tau)
&{\rm and}&\textstyle
\sum\limits_{\ell \in \m Z/d\m Z}\theta_{[1]}(\ell /d,\tau)=d\,\theta_{[1]}(0,d^2\tau).
\end{eqnarray*}
\el

\begin{proof}Just write the left-hand sides as a double sum over $m$ in $\m Z$
and $\ell $ in $\m Z/d\m Z$ and notice that 
$\sum_{\ell \in \m Z/d\m Z}e^{2i\pi \ell m/d}$ is equal to $d$ when $d$ divides $m$ 
and is equal to $0$ otherwise.
\end{proof}

The last formula deals with an element $\si
=\mbox{
\scriptsize 
	$\left(\!
	\begin{array}{cc} \al&\be   \\
	\ga&\de
	\end{array}\!
	\right)$} \in {\rm SL}(2,\m Z).
$

\bl 
\label{lemtrafor} {\bf Transformation formula} If $\si\equiv  \mathds{1}$ {\rm mod} $2$,
and $\ga>0$, then
\begin{eqnarray}
\label{eqntrafor}
\theta(0,\si\tau)
&=& i^\frac{\de-1}{2}(\!\tfrac{\ga}{\de}\!)\;
(\ga\tau+\de)^{\frac12}\;\theta(0,\tau).
\end{eqnarray}
\el
In this formula, the ${\rm SL}(2,\m Z)$ action 
on the upper half plane $\m H$ is the standard action 
$\si\tau =\frac{\al\tau+\be}{\ga\tau+\de}$,
the number $z^{\frac12}$
is the square root of a complex number $z\in \m H$ 
whose real part is non negative, and
the symbol  $(\!\frac{\ga}{\de}\!)=\pm 1$ is 
still the Jacobi symbol.

Note that Formula \eqref{eqntrafor} can be equivalently rewritten as
\begin{eqnarray*}
\theta(0,\si\tau)
&=& \eps_{\de}\,(\!\tfrac{2\ga}{\de}\!)\;
(\ga\tau+\de)^{\frac12}\;\theta(0,\tau),
\end{eqnarray*}
where $\eps_{\de}= 1$ when $\de\equiv 1$ mod $4$, 
and $\eps_{\de}=-i$ when $\de\equiv 3$ mod $4$.

\begin{proof} Up to sign, Formula \eqref{eqntrafor} follows from the 
following two formulas 
\begin{eqnarray*}
\theta(0,\tau+2)&=&\theta(0,\tau),\\
\theta(0,-1/\tau)&=&(-i\tau)^\frac{1}{2}\, \theta(0,\tau).
\end{eqnarray*}
and from the fact that the map $(\si,\tau)\mapsto \si\tau+\de$ 
is a cocycle on ${\rm SL}(2,\m Z)\times \m H$.
The precise determination of the sign is a classical issue
due to Hecke.
It can be found for instance in \cite[p.181]{Rad73} 
in \cite[p.148]{Kob84} or in \cite[p.32]{Mum83}.
\end{proof}

\br
\label{remmla}
This precise determination of the sign 
is important for us because it will allow us to decide 
whether the critical value we will find is 
$\la=\sqrt{a}+i\sqrt{b}$ or its opposite.
This is particularly important when $a$ is a square, because 
in this case $\la$ and $-\la$ are not Galois conjugate and one can not apply 
Proposition \ref{procrival}.i.
\er

The following corollary of Lemma \ref{lemtrafor} will be very useful.

\bc
\label{cortrafor}  
If $\si\equiv  \pm\mathds{1}$ {\rm mod} $4$,
 then, for all $\tau$ in $\m H$, one has
\begin{eqnarray}
\label{eqntrafor2}
\frac{\theta_{[0]}(0,\si\tau)}{\theta_{[0]}(0,\tau)}&=& \frac{\theta_{[1]}(0,\si\tau)}{\theta_{[1]}(0,\tau)}\, .
\end{eqnarray}
\ec

\begin{proof}
Let  
$$\si'
=\mbox{
$\left(\!
\begin{array}{cc} \al&\be'   \\
\ga'&\de
\end{array}\!
\right)$} 
\;\;{\rm  and}\;\; \si''
=\mbox{
$\left(\!
\begin{array}{cc} \al&\be''   \\
\ga''&\de
\end{array}\!
\right)$}, 
$$ 
with 
$\be'=2\be$, $\ga'=\ga/2$ and $\be''=\be/2$, $\ga''=2\ga$, 
so that 
$$
\si'(2\tau)=2\si\tau\;\;{\rm  and}\;\; \si''(\tau/2)=\tfrac12 \si\tau\, .
$$
Since the matrix $\si$ is equal to $\pm\mathds{1}$ mod $4$, the two matrices 
$\si'$ and $\si''$ are equal to  $\mathds{1}$ mod $2$. 
Therefore we can apply the transformation formula in Lemma \ref{lemtrafor}
to both pairs $(\si',2\tau)$ and $(\si'',\tau/2)$. 
Using the multiplicativity properties of the Jacobi symbol, we see that the 
following two ratios are given by the same formula
\begin{equation*}
\frac{\theta(0, 2\si\tau)}{\theta(0,2\tau)}=
\frac{\theta(0,\tfrac12 \si\tau)}{\theta(0,\tfrac12\tau)}.
\end{equation*}
We now conclude thanks to Equalities \eqref{eqnt0t1t2}.
\end{proof}

\subsection{The condition on theta contants} 

The first step in the proof of Proposition \ref{prorapirb} is 
the following criterion on $\la, \tau$ which ensures that the functions $f_{z,\tau}$ 
are $\la$-critical. This criterion is a  relation between  
``theta constants'', i.e. theta functions evaluated at $z=0$.

\bl
\label{lemrapirb}
Let $\tau\in \m H$ and $\la\in \m C$.
The function $\theta_\tau$ is $(\la,d)$-critical
if and only if 
one has the equalities
\begin{equation*}
\label{eqnrapirb}
\mbox{} \hspace{8em} 
\la= d\,\frac{\theta_{[0]}(0,d^2\tau)}{\theta_{[0]}(0,\tau)}= 
d\, \frac{\theta_{[1]}(0,d^2\tau)}{\theta_{[1]}(0,\tau)}\; .
\hspace{6em} (T_{\la,\tau})
\end{equation*}
\el

\begin{proof} For $w$ in $\m C$ we introduce the function
$$z\mapsto F_w(z)=F_w(z,\tau):=\theta(z+w,\tau)\,\theta(z-w,\tau).$$
We want to know when the two functions 
$\sum_\ell F_{\ell/d}$ and $F_0=\theta^2$ are proportional. 
The key point of the proof is that all these functions 
$F_w$ live in the same two-dimensional vector space 
and that this vector space has a very convenient basis:
$(\theta_{[0]},\theta_{[1]})$. 
We only have to express that the coefficients of our two functions in this basis are proportional. 
These coefficients are given by the following calculation in which we apply successively the addition formula and the isogeny formula,
\begin{eqnarray*}
\sum_\ell F_{\ell /d}(z,\tau )&=&
\sum_\ell  \theta_{[0]}(\ell /d,\tau)\; \theta_{[0]}(z,\tau) +
\sum_\ell  \theta_{[1]}(\ell /d,\tau)\; \theta_{[1]}(z,\tau)\\
&=& d\,\theta_{[0]}(0,d^2\tau)\; \theta_{[0]}(z,\tau) +
 d\, \theta_{[1]}(0,d^2\tau)\; \theta_{[1]}(z,\tau)
 \;\;\; {\rm and}
\end{eqnarray*}
\begin{eqnarray*}
\theta(z,\tau)^2&=&
\theta_{[0]}(0,\tau)\; \theta_{[0]}(z,\tau) +
 \theta_{[1]}(0,\tau)\; \theta_{[1]}(z,\tau).\hspace*{3em}
\end{eqnarray*}
These two functions are proportional with  proportionality factor $\la$
if and only if one has
\begin{equation*}
\label{eqnlat0t1} 
\la= d\,\frac{\theta_{[0]}(0,d^2\tau)}{\theta_{[0]}(0,\tau)}= 
d\, \frac{\theta_{[1]}(0,d^2\tau)}{\theta_{[1]}(0,\tau)}\; .
\end{equation*}
This is the criterion 
$(T_{\la,\tau})$.
\end{proof}

\subsection{The modular curve $X(4)$}
\label{secmodcur}

In order to interpret the condition $(T_{\la,\tau})$, the following classical description of the modular curve $X(4)$ will be very useful.

For $m\geq 1$, we introduce the principal congruence subgroup of level $m$
\begin{eqnarray*}
\label{eqnconsub} 
\Gamma(m)&:=&\{ \si\in {\rm SL}(2,\m Z)\mid \si\equiv \pm\mathds{1}\;\;{\rm mod}\; m \}
\end{eqnarray*}
and the modular curve of level $m$ 
\begin{eqnarray*}
\label{eqnmodcur} 
X(m)&:=&\Gamma(m)\backslash \m H .
\end{eqnarray*}
Note that  the element $-\mathds{1}\in {\rm SL}(2,\m Z)$ acts trivially on $\m H$.
It is classical that $X(m)$ is a Riemann surface with  finitely many cusps whose genus can be 
calculated thanks to Hurwitz formula. 
In this elementary paper we will only deal with $m=4$. In this case, 
$X(4)$ has genus zero and six cusps. The following lemma gives a nice interpretation of this fact. 

We introduce the meromorphic function $\Phi$ on $\m H$ given by, for all $\tau$ in $\m H$,
\begin{eqnarray*}
\label{eqnmapPhi} 
\Phi(\tau):=\frac{\th_{[1]}(0,\tau)}{\th_{[0]}(0,\tau)}.
\end{eqnarray*}

\bl 
\label{lemmodcur}
The map $\Phi$ induces a biholomorphism 
\begin{eqnarray*}
\label{eqnmapphi} 
\ph:X(4)&\longrightarrow &\m P^1\m C\smallsetminus\{ 0,\infty,\pm 1,\pm i\}.
\end{eqnarray*}
\el

This lemma tells us that, as an hyperbolic surface, $X(4)$ is the ''regular ideal octahedron''.

The statement of this lemma is equivalent to the following four facts 
on the meromorphic map $\Phi$.\\
$(a)$ For all $\si$ in $\Ga(4)$ and all $\tau$ in $\m H$ one has $\Phi(\si\tau) =\Phi(\tau)$.\\
$(b)$ For all $\tau$ in $\m H$, one has $\Phi(\tau)\neq 0,\infty,\pm 1,\pm i$.\\
$(c)$ If $\Phi(\tau)=\Phi(\tau')$, 
there exists $\si$ in $\Ga(4)$ such that $\tau'=\si\tau$.\\
$(d)$ For all $z\neq 0,\infty,\pm 1,\pm i$, there exists $\tau$ in $\m H$
with $\Phi(\tau)=z$.

Note that the first fact $(a)$ is the most important one for us in order to prove Theorem \ref{thmrapirb} and that 
it is just a restatement of Corollary \ref{cortrafor}.

\begin{proof}[Proof of Lemma \ref{lemmodcur}] This lemma is  classical for the experts. We will just relate it to the existing litterature.
We will deduce these four facts 
from a very similar statement in Mumford's book \cite{Mum83}. 
In this book, Mumford uses the four Jacobi theta-functions $\theta_{a,b}$ with $a$, $b$ equal to $0$ or $1$, given by
$$\textstyle
\theta_{a,b}(z,\tau)=\sum\limits_{m\in \m Z} e^{i\pi\tau(m+\frac{a}{2})^2}e^{2i\pi(m+\frac{a}{2})(z+\frac{b}{2})}
$$
It is proven in 
\cite[Theorem 10.1 p.51]{Mum83} that the map $\Psi$ given in homogeneous coordinates by
$$
\Psi:\tau\mapsto [\th_{0,0}^2(0,\tau),\th_{0,1}^2(0,\tau),\th_{1,0}^2(0,\tau)]
$$
induces an biholomophism $\psi$ between  the curve $X(4)$ and the 
curve
$$
C:=\{[x_0,x_1,x_2]\in \m P^2\m C\mid x_0^2=x_1^2+x_2^2\; \mbox{\rm and all}\; x_i\neq 0\}
$$ 
which is a conic with six points removed.
By the addition formula \eqref{eqnaddfor}, these theta-constants $\theta_{a,b}$  are related to the theta-constants 
$\theta_{[0]}$ and 
$\theta_{[1]}$\, :
\begin{eqnarray*} 
\theta_{0,0}^2(0,\tau)&=&\theta_{[0]}^2(0,\tau)+\theta_{[1]}^2(0,\tau),\\
\theta_{0,1}^2(0,\tau)&=&\theta_{[0]}^2(0,\tau)-\theta_{[1]}^2(0,\tau) ,\\
\theta_{1,0}^2(0,\tau)&=&2\,\theta_{[0]}(0,\tau)\,\theta_{[1]}(0,\tau).
\end{eqnarray*} 
Hence one can express in a simple way the map $\Psi$ thanks to the function $\Phi$:
$$ \Psi(\tau)=[1\!+\!\Phi^2(\tau)\, ,\, 1\!-\!\Phi^2(\tau)\, ,\, 2\Phi(\tau)],
$$
for all $\tau$ in $\m H$. 
\end{proof}

\br
\label{remmodcur} Note that this identification  of 
$X(4)$ is equivariant. More precisely,
the finite group $G:={\rm PGL}(2,\m Z/4\m Z)$ 
has cardinality $48$ and acts by biholomorphisms or biantiholomorphisms on $X(4)$.
The biholomorphism $\ph$ identifies this group $G$ with the group of isometries of the octahedron.
This follows from the identities
$$
\Phi(-\ol{\tau})=\ol{\Phi(\tau)}
\;\; ,\;\;
\Phi(\tau+1)=\Phi(\tau)
\;{\rm and}\;\;
\Phi(-1/\tau)=
\frac{-\Phi(\tau)+1}{\Phi(\tau+1}\, .
$$
\er

\subsection{Elliptic curves with complex multiplication} 

We can now go on the proof of Proposition \ref{prorapirb}, 
by explaining how we will check that a pair $(\la,\tau)$ satisfies
Condition $(T_{\la,\tau})$.

For $\tau$ in $\m H$, we introduce the lattice $\La_\tau=\m Z\tau\oplus\m Z 1$ of $\m C$
so that the quotient $E_\tau:=\m C/\La_\tau$ 
is the elliptic curve associated to $\tau$.
We will see that 
the values of $\la$ and $\tau=\tau_{k,p}$ in Theorem \ref{thmrapirb} 
have been chosen so that
the elliptic curve $E_{\tau}$ has complex multiplication by 
$\mu:=\overline{\la}^2$.
See \cite{Sche10} for more classical applications of complex multiplication.
More precisely, they have been chosen so that $\mu\La_{\tau}=\La_{d^2\tau}$. 
This means that
\begin{eqnarray}
\label{eqncommul1}
d^2\tau&=& \mu \; (\al\tau+\be),\\
\label{eqncommul2}
1&=& \mu\;(\ga\tau+\de),
\end{eqnarray}
for a matrix 
$\si
=\mbox{
$\left(\!
\begin{array}{cc} \al&\be   \\
\ga&\de
\end{array}\!
\right)$} \in {\rm SL}(2,\m Z)
$. We will be able to impose on $\si$  the extra condition  
\begin{eqnarray} 
\label{eqng0mod4}
\ga>0
&{\rm and}& \si\equiv \pm\mathds{1}\;\;{\rm mod}\; 4.
\end{eqnarray}
We explain why such a choice is relevant in the following lemma.

\bl
\label{lemmulcom}
Let $\tau\in \m H$ and $d$ an odd integer.\\ 
$(a)$ The function $\theta_\tau$ is $(\la,d)$-critical 
for some $\la$ in $\m C$ if and only if there exists 
$\si
=\mbox{
\scriptsize 
$\left(\!
\begin{array}{cc} \al&\be   \\
\ga&\de
\end{array}\!
\right)$} \in {\rm SL}(2,\m Z)
$ satisfying \eqref{eqng0mod4} such that $g\tau=d^2\tau$.\\
$(b)$
In this case, 
the critical value  is given by, setting $\mu=(\ga\tau+\de)^{-1}$,
\begin{equation}
\label{eqnlamuba}
\la=\eps_{\de}\,(\!\tfrac{\ga}{\de}\!)\;\ol{\mu}^{1/2}.
\end{equation} 
\el

\noindent
Recall that  $\eps_{\de}= 1$ when $\de\equiv 1$ mod $4$, 
and $\eps_{\de}=-i$ when $\de\equiv 3$ mod $4$.

\begin{proof}[Proof of Lemma \ref{lemmulcom}] 
$(a)$ According to Lemma \ref{lemrapirb},
the function $\theta_\tau$ is $(\la,d)$-critical
for some $\la$ in $\m C$ if and only if 
one has the equalities $\Phi(\tau)=\Phi(d^2\tau)$.

According to Lemma \ref{lemmodcur}, this equality 
is equivalent to the existence of an element
$\si$ in $\Ga(4)$ such that $\si\tau=d^2\tau$.
Note that this forces the lower-left coefficient 
$\ga$ to be non zero.
Replacing $\si$ by $-\si$ if necessary, we can also assume that $\ga >0$.

$(b)$ To compute the critical value $\la$, 
we use again Lemma \ref{lemrapirb} combined with Formulas \eqref{eqnt0t1t2}. We obtain
\begin{equation*}
\la=d\,\frac{\theta_{[0]}(0,d^2 \tau)}{\theta_{[0]}(0,\tau)}=d\,
\frac{\theta(0,d^2 2\tau)}{\theta(0,2\tau)}
\; . 
\end{equation*}
Let  
$\si'
=\mbox{
$\left(\!
\begin{array}{cc} \al&\be'  \\
\ga'&\de
\end{array}\!
\right)$},
$ 
with 
$\be'=2\be$, $\ga'=\ga/2$, 
so that 
$
\si'(2\tau)=d^2 2\tau.
$
Since the matrix $\si$ equals  $\pm\mathds{1}$ mod $4$, the matrix 
$\si'$ is equal to  $\mathds{1}$ mod $2$.
Hence 
we can apply the transformation formula in Lemma \ref{lemtrafor}
to  $\si'$ and $\tau'$. 
We get
\begin{equation}
\label{eqnlaepde}
\la =d\, \eps_{\de}\,(\!\tfrac{\ga}{\de}\!)\;(\ga\tau+\de)^\frac12.
\end{equation}

By assumption, the number $\mu:=(\ga\tau+\de)^{-1}$ is a complex multiplication for the lattice $\La_\tau$,
 more precisely, one has
$\mu\La_\tau=\La_{d^2\tau}$. This gives
the equality $\mu\ol{\mu}=d^2$, and Equation
\eqref{eqncommul2} can be rewritten as 
$$
d^2\, (\ga\tau+\de)=\ol{\mu}.
$$ 
Now  Equation \eqref{eqnlaepde} 
can be rewritten as
$\la=\eps_{\de}\,(\!\frac{\ga}{\de}\!)\,\ol{\mu}^\frac12$.
\end{proof}

\subsection{Choosing the elliptic curve} 
\label{secchoell}

In order to end the proof of Proposition 
\ref{prorapirb},
it only remains to characterize those points $\tau$ in $\m H$ that satisfy the condition
$\si\tau=d^2\tau$ for some $\si$ in $\Ga(4)$ and to express the $d$-critical value $\la$
given in Lemma \ref{lemmulcom} without using $\si$. 
\vs 

We begin by recalling the notation of Proposition \ref{prorapirb}.
Let $a$, $b$ be positive integers with 
$a\!\equiv\!\frac{(d+1)^2}{4}\;{\rm  mod}\; 4$
and $a\!+\!b\!=\!d$.
We introduced the ``fundamental parameter'' 
$\tau_0:=\tfrac{1}{4d^2}(a-b-d^2+2i\sqrt{ab})\in \m H$. 
For $k$ in $\m Z$, we introduced the  integer $N_k:=d^2|k+\tau_0|^2$.
For $p\in \m N$ dividing $N_k$, we also introduced the ``associated parameters''
$\tau_{k,p}:=(k+\tau_0)/p\in \m H$. 

\bl
\label{lemtaukg0}
(a) Let $\tau$ in $\m H$ such that there exists $\si$ in $\Ga(4)$ 
for which $\si\tau=d^2\tau$,
then there exists $k, p$ as above such that $\tau=\tau_{k,p}$ or 
$\tau=-\ol{\tau}_{k,p}$.\\
$(b)$ Conversely, for every $\tau=\tau_{k,p}$ as above, there exists $\si$ in $\Ga(4)$
such that  $\si\tau=d^2\tau$.\\
$(c)$ This matrix $\si$
can be chosen to be 
$\si=\mbox{
\scriptsize 
$\left(\!
\begin{array}{cc} 2(a\!-\!b)\!-\!d^2(1\!-\!4k)&-4 N_k/p   \\
4p&1-4k
\end{array}\!
\right)$} .$
\el

\begin{proof} We first notice that, in this lemma, we can always add the extra conditions $\ga>0$ 
and $\si\equiv \mathds{1}$ mod $4$. Indeed, if needed, we can always replace
the matrix  $\si$ by $-\si$ without changing 
the point $\tau$. We can also replace $\si$ by the matrix 
 $\mbox{
\scriptsize 
$\left(\!
\begin{array}{cc} \al&-\be   \\
-\ga&\de
\end{array}\!
\right)$} $, the point $\tau$ is then replaced by $-\ol{\tau}$.
\vs 
 
$(a)$ We set $\mu:=(\ga\tau\!+\!\de)^{-1}$.
Since the matrix $\si=\mbox{
\scriptsize 
$\left(\!
\begin{array}{cc} \al&\be   \\
\ga&\de
\end{array}\!
\right)$} $ has deter\-minant $1$,
the equations \eqref{eqncommul1} and \eqref{eqncommul2} 
can be rewritten as 
\begin{eqnarray}
\label{eqncommul3}
\;\; \ga^{-1}(\mu^{-1} -\de)&=& \tau,\\
\label{eqncommul4}
\mu^2- t_0\mu+d^2&=& 0,\;\; 
\end{eqnarray}
with $t_0:=\al+d^2\de$. We introduce
$a$, $b$ such that $t_0=2(a-b)$ 
and $d=a+b$.\\
$\star$ Since $\al\equiv \de\equiv  1$ mod $4$ and $d^2\equiv 1$ mod $4$, 
both numbers  $a$ and $b$ are integers.\\ 
$\star$ Since $\mu$ is not a real number, one has $|t_0|<2d$ and the integers $a$ and $b$ are positive.\\ 
$\star$ Since $\si\equiv \mathds{1}$ mod $4$, one has $\al\de\equiv 1$ mod $16$, and these integers satisfy 
$a-b\equiv \frac{1+d^2}{2}$ mod $8$,
and hence 
$a\!\equiv\!\frac{(d+1)^2}{4}\;{\rm  mod}\; 4$.\\
$\star$ Since $Im(\mu)<0$, solving Equation \eqref{eqncommul4}, one gets the equality 
\begin{eqnarray}
\label{eqnmu0la0} 
\mu&:=\;(\sqrt{a}-i\sqrt{b})^2
\; =\; a-b-2i\sqrt{ab}
\end{eqnarray}
$\star$ We write $\de=1-4k$ with $k\in \m Z$, and one computes
\begin{eqnarray*}
N_k=d^2|k+\tau_0|^2&=&
\tfrac{1}{16d^2}[((a-b)-(1-4k)d^2)^2+4ab]\\
&=&\tfrac{1}{16}[1-2(a-b)(1-4k)+d^2(1-4k)^2]\\
&=&(1-\al\de)/16 \; =\; -\be\ga/16.
\end{eqnarray*}
$\star$ We write  $\ga=4p$ so that this integer $p$
is a divisor of the integer $N_k$.\\
$\star$ Plugging these informations in \eqref{eqncommul3} gives 
$\tau=(k+\tau_0)/p$.
\vs 

$(b)$  and $(c)$\; We assume now that $\tau=\tau_{k,p}$ as above and 
we want to construct the matrix $\si$. 
We follow the same computation as in $(a)$ in opposite order.
We set $\mu:=a-b-2i\sqrt{ab}$ and
$t_0:=2(a-b)$
so that Equation \eqref{eqncommul4} is satisfied.
We choose $\de:=1-4k$ and 
$\al:=t_0-d^2\de$. We first note that 
\begin{equation}
\label{eqnadmod4}
\al\de\equiv 1\;\;{\rm mod}\;\; 16.
\end{equation}
To check \eqref{eqnadmod4}, just remember that one has 
$t_0\equiv 1+d^2$ mod $16$, and hence
$$
\al\de-1\equiv (\de-1)(1-\de d^2)\equiv 4^2\equiv 0\;\;\;{\rm mod}\; 16.
$$
Then the same computation as above gives the equality $(1-\al\de)/16=N_k$.
Hence if we choose   $\ga:=4p$ and $\be:=-4N_k/p$,
the matrix $\si$ is in $\Ga(4)$.
By construction these coefficients satisfy also Equality  \eqref{eqncommul3}.
Hence the matrix $\si$ satisfies Equalities \eqref{eqncommul1} and \eqref{eqncommul2}.
\end{proof}

\begin{proof}[Proof of Proposition \ref{prorapirb}] 
This proposition is now just a straightforward consequence of 
Lemmas \ref{lemmulcom} and \ref{lemtaukg0} combined with Formula
\eqref{eqnmu0la0}.
\end{proof}

\subsection{Conclusion and Perspective} 
\label{secconclusion}

The aim of this paper was to explain why 
the algebraic integers $\sqrt{a}+i\sqrt{b}$ that occur 
in the lists of section \ref{secnumexp} are indeed 
$d$-critical values when $a,b$ are positive integers
with $a+b=d$ and $a\!\equiv\frac{(d+1)^2}{4}$ mod $4$.
To keep this paper as elementary as possible, we have only discussed here these $d$-critical values.
\vs 

However, in the lists of section \ref{secnumexp}, there are still
remaining intriguing $d$-critical values. 
In a more technical forthcoming paper \cite{CSAGII}, 
we will see that these  $d$-critical values belong to a 
wide class of critical values on finite abelian groups 
that can be explained by an extension of the construction of Proposition \ref{prorapirb}. 
This will be a nice application of the 
abelian varieties and their theta funtions, 
relying on works of Siegel, Stark, Igusa and of Taniyama-Shimura. See \cite{Beau13} 
for a recent paper surveying previous applications of these tools.
Indeed we will prove in \cite{CSAGII}.

\bt
\label{thmmaivap}
Let $A=\m C^g/\La$ be a principally polarized abelian variety
and  $\nu$ be a unitary $\m Q$-endomorphism of $A$ preserving a theta structure of level $2$. Let $T_\nu$ be its tangent map, $G_\nu$  the group $\La/(\La\cap T_\nu\La)$
and $d_\nu$  the order of $G_\nu$. 
Then there is a critical value 
$\la_\nu=\ka_\nu\, d_\nu^{1/2}\,{\rm det}_{\m C}(T_\nu)^{1/2}$
on $G_\nu$ with $\ka_\nu^4=1$.
\et

\noindent
Note that $|\la_\nu|=d_\nu^{1/2}$ and that one can compute the fourth root of unity $\ka_\nu$.

Note also that, any finite abelian group can occur as a group $G_\nu$ but,
even when $g>1$, this group $G_\nu$ may be cyclic. 
We will construct concrete examples with $G_\nu$ cyclic 
when the abelian variety $A$ has complex multiplication.
These constructions will explain all the intriguing $d$-critical values in our lists.
For instance we will show in \cite{CSAGII}
by using abelian surfaces.

\bc
\label{corparacri}
Let $d_j\!=\!a_j\!+\!b_j$  with $d_1\wedge d_2=1$ and
$a_j\!\equiv\! \frac{(d_j+1)^2}{4}\!+\! 2$ {\rm mod}~$4$ be positive integers.
Then  
$\la\!=\!(\sqrt{a_1}\!+\! i\sqrt{b_1})(\sqrt{a_2}\!-\! i\sqrt{b_2})
$ 
is $d$-critical.
\ec

\bc
\label{corrarbrc} 
Let $d=a\!+\! b\!+\!c$ be positive integers, with  $b^2>4ac$, 
and  $a\!\equiv\! b\!\equiv\! c\!\equiv\!   1$  mod $4$.
Then 
$
\la= \sqrt{a}\!+\!\sqrt{c} + i\sqrt{b\!-\! 2 \sqrt{ac}} 
$ 
is  $d$-critical.
\ec
\noindent 
A key remark in the proof of  Corollary \ref{corrarbrc}, is an old factorization formula:
\begin{equation*} 
\label{eqn1pr5pi}
\textstyle
\sqrt{a}\!+\!\sqrt{c}\!+\!i\sqrt{b\!-\!2\sqrt{ac}}
\; =\;
\sqrt{a}\left(\!1\!+\! i\sqrt{\tfrac{b+\sqrt{b^2-4ac}}{2a}}\right)
\left(\!1\!-\! i\sqrt{\tfrac{b-\sqrt{b^2-4ac}}{2a}}\right).
\end{equation*}

{\small
\bibliography{theta}
}
\vs 

{\small
\noindent
Y. \textsc{Benoist}: CNRS, 
Universit\'e Paris-Saclay,\hfill
\texttt{yves.benoist@u-psud.fr}}

\end{document}